\documentclass[12pt]{article}
\usepackage{amsmath}
\usepackage{color}
\usepackage{fancyhdr}
\usepackage{verbatim}
\usepackage{amsfonts,longtable,mathrsfs}
\usepackage{epsfig}
\usepackage{amsfonts}
\usepackage{color}
\usepackage[numbers,sort&compress]{natbib}

\newtheorem{definition}{Definition}
\begin{document}
\begin{center} \textbf{ A group theory approach towards some rational difference equations
  }\end{center}
\medskip
\begin{center}
Mensah Folly-Gbetoula\footnote{Mensah.Folly-Gbetoula@wits.ac.za}, Nkosingiphile  Mnguni and  A. H. Kara\vspace{1cm}
\\ School of Mathematics, University of the Witwatersrand, Private Bag 3, Johannesburg 2050, South Africa.
\end{center}
\begin{abstract}
\noindent
 A full Lie point symmetry analysis of rational difference equations is performed. Non-trivial symmetries are derived and exact solutions using these symmetries are obtained.
\end{abstract}
\section{Introduction} \setcounter{equation}{0}
Over a century ago, symmetries became a centre of interest of several authors after the work of Sophus Lie \cite{Lie} on differential equations. He studied the continuous group of transformations that leaves the differential equations invariant. This concept  of symmetries is strongly related to the existence of conservation laws and the relationship between them has attracted great interest among researchers following the work of Noether \cite{Noether}. The extension of this idea to difference equations is now well-documented (see \cite{Hydon} and references herein). In \cite{Hydon}, Hydon developed a symmetry based algorithm enabling one to derive solutions of difference equations without making any special lucky guesses. Hydon  emphasized on second-order difference equations, although his algorithm is valid for any order. When it comes to higher-order equations, computations can be cumbersome and extra ansatz may be needed to ease the calculations.
\par \noindent We aim to extend the work by Elsayed \cite{Elsayed} where the author investigated the dynamics and solutions of
\begin{align}\label{1.0}
x_{n+1}=\frac{x_{n-5}}{\pm 1\pm x_{n-1}x_{n-3}x_{n-5}},
\end{align}
where the initial conditions $x_{-5},\, x_{-4},\,x_{-3},\,x_{-2},\,x_{-1}$ and $x_0$ are arbitrary non-zero real numbers. For related work, see \cite{Mensah, FK}. One can notice that equations (\ref{1.0}) are just special cases of a more general form
\begin{align}\label{1.1}
x_{n+1}=\frac{x_{n-5}}{a_n +b_n x_{n-1}x_{n-3}x_{n-5}},
\end{align}
where $(a_n)_{n\in \mathbb{N}_0}$ and $(b_n)_{n\in \mathbb{N}_0}$ are arbitrary sequences. We will use a symmetry based method to solve (\ref{1.1}). Equivalently, we study
\begin{align}\label{1.2}
u_{n+6}=\frac{u_n}{A _n + B _n u_nu_{n+2}u_{n+4}}
\end{align}
instead, where $(A_n)_{n\in \mathbb{N}_0}$ and $(B_n)_{n\in \mathbb{N}_0}$ are arbitrary sequences. This means we can only compare $x_i$ with $u_{i+5}$. Furthermore, we use the same technique to obtain exact solutions of
\begin{align}\label{xn5}
x_{n+1}=\frac{x_{n-4}x_{n-3}x_{n-2}}{x_{n-1}x_{n}(\lambda + \mu x_{n-4}x_{n-3}x_{n-2})},
\end{align}
where $\lambda, \mu \in \mathbb{R}$, and again we study
\begin{align}\label{un5'}
u_{n+5}=\frac{u_nu_{n+1}u_{n+2}}{u_{n+3}u_{n+4}(\lambda + \mu u_nu_{n+1}u_{n+2})},
\end{align}
instead. Note that solutions of \eqref{xn5} were found in \cite{stevo2}; however, their method is completely different from ours.
\section{Definitions and notation}
The definitions are taken from Hydon \cite{Hydon} and most of the notation follows from the same book.
\begin{definition}
A parameterized set of point transformations,
\begin{equation}
\Gamma_{\varepsilon} :x\mapsto \hat{x}(x;\varepsilon),
\label{eq: b}
\end{equation}
where $x=x_i, $ $i=1,\dots,p$ are continuous variables, is a one-parameter local Lie group of transformations if the following conditions are satisfied:
\begin{enumerate}
\item $\Gamma_0$ is the identity map if $\hat{x}=x$ when $\varepsilon=0$
\item $\Gamma_a\Gamma_b=\Gamma_{a+b}$ for every $a$ and $b$ sufficiently close to 0
\item Each $\hat{x_i}$ can be represented as a Taylor series (in a neighborhood of $\varepsilon=0$ that is determined by $x$), and therefore
\end{enumerate}
\begin{equation}
\hat{x_i}(x; \varepsilon)=x_i+\varepsilon \xi _i(x)+O(\varepsilon ^2), i=1,...,p.
\label{eq: c}
\end{equation}
\end{definition}
 Consider the $p$th-order difference equation
\begin{equation}
u_{n+p}=\Omega(n, u_n, \dots, u\;_{n+p-1}),
\label{eq: a}
\end{equation}
for some function $\Omega$. Assume  the point transformations are of the form
\begin{equation}
\hat{n}=n ;\qquad
\hat{u}_n = u_n + \varepsilon Q(n, u_n)+O(\varepsilon ^2)
\label{eq: d}
\end{equation}
with the corresponding infinitesimal symmetry generator
\begin{align}
\begin{split}
X=&Q(n,u_{n})\frac{\partial}{\partial u_{n}}+SQ(n,u_{n})\frac{\partial}{\partial u_{n+1}}+\dots
+S^{(p-1)}Q(n,u_{n})\frac{\partial}{\partial u_{n+p-1}},
\label{eq: e}
\end{split}
\end{align}
where $S$ is  the shift operator, i.e., $S:n\mapsto n+1$.
The symmetry condition is defined as
\begin{equation}
\hat{u}_{n+p}=\Omega(n, \hat{u}_{n}, \hat{u}_{n+1},...,\hat{u}_{n+p-1}),
\label{eq: g}
\end{equation}
whenever \eqref{eq: a} is true. Substituting the Lie point symmetries \eqref{eq: d} into the symmetry condition \eqref{eq: g} leads to  the linearized symmetry condition
\begin{equation}\label{LSC}
S^{(p)}Q-X\Omega=0
\end{equation}
whenever \eqref{eq: a} holds.
\begin{definition}
$V_n$ is invariant under the Lie group of transformations \eqref{eq: d} if and only if $X V_n =0$.
\end{definition}
We define the functions
\begin{align}\label{Theta}
\Theta ^n (\theta_j,s,l)=\prod_{j=s}^{l}\theta_{2j+n},\quad n=0,1,
\end{align}
and we adopt the standard conventions
\begin{align}
\prod_{j=s}^l \theta_j=1 \text{ when } s>l \text{ and } \sum _{j=n}^{l_0}\theta_j=0  \text{ when } n>l_0.
\end{align}
We refer the reader to \cite{Hydon} for a deeper understanding of the concept of symmetry analysis of difference equations.
\section{Main results}
\subsection{On the difference equations (\ref{1.2})}
Consider the sixth-order difference equations of the form \eqref{1.2}, that is,
\begin{align}\label{1.3}
u_{n+6}=\frac{u_n}{A _n + B _n u_nu_{n+2}u_{n+4}}.
\end{align}
We impose the symmetry condition (\ref{LSC}) and simplify the resulting equation to get
\begin{align}\label{LSCa}
&S^6Q+
\frac{B_n{u_n}^2u_{n+2}S^4Q +B_n{u_n}^2{u_{n+4}}S^2Q-{A_n} Q   } {( {B_n}{u_n} u_{n+2}{u_{n+4}} +A_n)^2}=0.
\end{align}
To solve for $Q$, we first differentiate  \eqref{LSCa} with respect to $u_n$ ( keeping $\Omega$ fixed and viewing $u_{n+2}$ as a function of $u_n, u_{n+4}$ and $\Omega$). This leads, after simplification, to
\begin{align}\label{A2}
\begin{split}
&S^4Q+ u_{n+4}(S^2Q)' - u_{n+4}  Q'+ \frac{2 u_{n+4} }{u_{n}}Q=0,
\end{split}
\end{align}
where $'$ denotes the derivative with respect to the independent variable.
We then differentiate (\ref{A2}) with respect to $u_{n}$ to get
\begin{align}\label{a7}
\begin{split}
 & Q''\left(n, u_n\right)- \frac{2  }{u_n}Q'\left(n, u_n\right) + \frac{2 }{u_n^{2}}Q\left(n, u_n\right)=0.
 \end{split}
\end{align}
The  solutions of \eqref{a7} are given by
\begin{align}\label{aa6'}
\begin{split}
Q\left(n,{u_{n}}\right) = \alpha_n  {u_n}^2 +\beta_n  {u_n}
\end{split}
\end{align}
for some functions $\alpha_n$ and $\beta_n$ of $n$. To obtain more information on $\alpha_n$ and $\beta_n$, we substitute (\ref{aa6'}) in (\ref{LSCa}) and we split the resulting equation to get the following:
\begin{subequations}\label{a12'}
\begin{align}\label{a12}
& {u_n}{u_{n+2}}{u_{n+4}}^2&:&\;  B_n\alpha_{n+4}=0\\
&{u_n}{u_{n+2}}^2{u_{n+4}}&:& \;B_n\alpha_{n+2}=0\\
&{u_n}{u_{n+2}}{u_{n+4}}&:& \; B_n(\beta_{n+2}+\beta_{n+4}+\beta_{n+6}) =0\\
&{u_{n}}&:& \;-A_n\alpha_{n}+\alpha_{n+6}=0  \\
&1&:&\;A_n\left( \beta_{n}-\beta_{n+6}\right) =0.
\end{align}
\end{subequations}
These equations (\ref{a12'}) reduce to
\begin{align}
&\alpha _n =0,\; \beta_{n+4}+\beta_{n+2}+\beta_n=0\label{ab}.
\end{align}
The expression of $\beta_n$ in \eqref{ab} is merely obtained by solving the corresponding characteristic equation $r^4+r^2+1=0$ for $r$. Assuming that $r_i,\; i=1,2,3,4$, are the solutions of this characteristic equation, then  $\beta_n$ is a linear combination of the $r_i^n$'s. In other words, the solutions of \eqref{ab}
are
\begin{align}
\alpha _n =0 \qquad \text{ and } \qquad \beta _n=(-1)^{n}{\beta}^{n}c_1
+(-1)^{n}{\bar{\beta}}^{n}c_2+{\bar{\beta}}^{n}c_3+{\beta}^{n}c_4
\end{align}
for some arbitrary constants $c_i,\; i=1,\dots,4$, and where $\beta= \exp{\left({\pi i}/{3}\right)}$. Thus, we obtain four characteristics given by
\begin{eqnarray}\label{generr}
Q_1=(-1)^{n}{\beta}^{n} u_n,\quad Q_2=(-1)^{n}{\bar{\beta}}^{n} u_n \quad Q_3={\bar{\beta}}^{n} u_n,\quad Q_4={\beta}^{n} u_n.
\end{eqnarray}
The four corresponding symmetry  generators $X_1,\,X_2,\,X_3$ and $X_4$ are given by
\begin{subequations}
\begin{eqnarray}\label{gener'}
\begin{split}
X_1=&  (-1)^{n}{\beta}^{n} u_n \partial u_n-  (-1)^{n}{\beta}^{n+1} u_{n+1} \partial_{u_{ n+1}}+(-1)^{n}{\beta}^{n+2} u_{n+2} \partial_{u_{ n+2}}-\\& (-1)^{n}{\beta}^{n+3} u_{n+3} \partial_{u_{ n+3}}+ (-1)^{n}{\beta}^{n+4} u_{n+4} \partial_{u_{ n+4}}- (-1)^{n}{\beta}^{n+5} u_{n+5} \partial_{u_{ n+5}},\label{ggener1}
\end{split}
\\
\begin{split}
X_2=&   (-1)^{n}{\bar{\beta}}^{n} u_n \partial u_n-  (-1)^{n}{\bar{\beta}}^{n+1} u_{n+1} \partial_{u_{ n+1}}+(-1)^{n}{\bar{\beta}}^{n+2} u_{n+2} \partial_{u_{ n+2}}-\\& (-1)^{n}{\bar{\beta}}^{n+3} u_{n+3} \partial_{u_{ n+3}}+ (-1)^{n}{\bar{\beta}}^{n+4} u_{n+4} \partial_{u_{ n+4}}- (-1)^{n}{\bar{\beta}}^{n+5} u_{n+5} \partial_{u_{ n+5}}\label{ggener2},
\end{split}\\
\begin{split}
X_3=&{\bar{\beta}}^{n} u_n \partial u_n+{\bar{\beta}}^{n+1} u_{n+1} \partial_{u_{ n+1}}+{\bar{\beta}}^{n+2} u_{n+2} \partial_{u_{ n+2}}+ {\bar{\beta}}^{n+3} u_{n+3} \partial_{u_{ n+3}}+\qquad  \\& {\bar{\beta}}^{n+4} u_{n+4} \partial_{u_{ n+4}}+{\bar{\beta}}^{n+5} u_{n+5} \partial_{u_{ n+5}},\label{gener3}
\end{split}
\\
\begin{split}
X_4=& {\beta}^{n} u_n \partial u_n+  {\beta}^{n+1} u_{n+1} \partial_{u_{ n+1}}+ {\beta}^{n+2} u_{n+2} \partial_{u_{ n+2}}+  {\beta}^{n+3} u_{n+3} \partial_{u_{ n+3}}+\qquad\\&  {\beta}^{n+4} u_{n+4} \partial_{u_{ n+4}}+  {\beta}^{n+5} u_{n+5} \partial_{u_{ n+5}}.\label{gener4}
\end{split}
\end{eqnarray}
\end{subequations}
Here, using $Q_4$, we introduce the canonical coordinate \cite{JV}
\begin{align}\label{sn}
s_n=\int{\dfrac{du_n}{{{\beta}^n}u_n}}=\frac{1}{\beta^n}\ln |u_n|.
\end{align}
In view of \eqref{ab}, we introduce the variable
\begin{align}\label{rn}
r_n=\beta ^{n+4}s_{n+4}+\beta ^{n+2}s_{n+2}+\beta ^n s_n.
\end{align}
It is easy to check that
\begin{align}
X_i\, r_n =\beta ^{n+4} +\beta ^{n+2} + \beta ^n =0, \quad i=1, \dots, 4.
\end{align}
Therefore $r_n$ is  invariant under $X_4$.
It is advantageous to use
\begin{align}\label{rn'}
|\tilde{r}_n|=\exp\left(-r_n\right),
\end{align}
that is, $ \tilde{r}_n=\pm 1/(u_nu_{n+2}u_{n+4})$. Here, we choose to use the plus sign and we have shown, using \eqref{1.2}, that
\begin{align}\label{rn''}
\tilde{r}_{n+2}=A_n \tilde{r}_{n}+B_n.
\end{align}
Therefore,
\begin{align}\label{rn'''}
\tilde{r}_{2n+k}=\tilde{r}_k\left( \prod_{k_1=0}^{n-1}A_{2k_1+k}\right)+\sum _{l=0}^{n-1}\left(   B_{2l+k}\prod _{k_2=l+1}^{n-1}A_{2k_2+k}\right), \quad k=0,1.
\end{align}
It is worthwhile to mention that equations in \eqref{rn'''} give the solution of  \eqref{rn''} for all $n$. By reversing all the change of variables, we have
\begin{align}
|u_n| 
=&\exp\Bigg[\left(-\beta\right)^n c_5 + \left(\beta\right)^n c_6 +
\left(\bar{\beta}\right)^n c_7 + \left(-\bar{\beta}\right)^n c_8 +\nonumber\\
&
 \left(-\beta\right)^{n+1} \sum\limits_{k_1 =0}^{n - 1}\frac{-i\sqrt{3}}{6} \left(-\bar{\beta}\right)^{k_1} r_{k_1}
 + \left(\beta\right)^{n+1} \sum_{k_2 =0}^{n - 1} \frac{i\sqrt{3}}{6}\left(\bar{\beta}\right)^{k_2} r_{k_2} +\nonumber \\ &
  (\bar{\beta})^{n+1} \sum_{k_3 =0}^{n - 1} \frac{-i\sqrt{3}}{6}\beta ^{k_3} r_{k_3} +
   (-\bar{\beta})^{n+1} \sum_{k_4 =0}^{n - 1} \frac{i\sqrt{3}}{6}(-\beta )^{k_4}r_{k_4}                      \Bigg]\nonumber\\
 =&{\Gamma _n}\exp\Bigg(\sum\limits_{k =0}^{n - 1}\frac{\sqrt{3}}{3}\left[(-1)^{n+k}+1\right]\text{Im}\left[ \gamma(n+1,k)\right]\ln |\tilde{r}_{k}|\Bigg),\label{solun}
\end{align}
where $\tilde{r}_k$ is given in \eqref{rn'''}, $\Gamma_n=\exp\{\left(-\beta\right)^n c_5 + \beta ^n c_6 +\bar{\beta} ^n c_7 + \left(-\bar{\beta}\right)^n c_8\}$ and $\gamma(l,k)=\beta ^l \bar{\beta}^{k}$.
\\
%\begin{rem}
\textbf{Note}. Equation (\ref{solun}) gives the solution of (\ref{1.2}) in a unified manner.\par \noindent
We can simplify \eqref{solun} further by splitting it into six categories. We have
\begin{align}\label{u6n}
 |u_{6n}| =&{\Gamma _{6n}}\exp\Bigg(\sum\limits_{k =0}^{6n - 1}\frac{\sqrt{3}}{3}\left[(-1)^{6n+k}+1\right]\text{Im}\left[ \gamma(6n+1,k)\right]\ln |\tilde{r}_{k}|\Bigg)\\
u_{6n}=&u_0 \prod_{s=0}^{n-1}\frac{\tilde{r}_{6s}}{\tilde{r}_{6s+2}}.
\end{align}
Similarly, after a straightforward but lengthy computation, we get
\begin{align}\label{8c}
&u_{6n+i}=u_i \prod_{s=0}^{n-1}\frac{\tilde{r}_{6s+i}}{\tilde{r}_{6s+2+i}}, i=0, \dots, 5.
\end{align}
\textbf{Note}. We can obtain (\ref{8c}) using (\ref{rn'}) which need not the use of absolute values. \par \noindent
Now, using \eqref{rn'''} and \eqref{8c}, we obtain the solutions of \eqref{1.2} as follows:
\begin{subequations}
\begin{align}\label{genesolun}
u_{6n}=&u_0 \prod_{s=0}^{n-1}\frac{\left( \prod\limits_{k_1=0}^{3s-1}A_{2k_1}\right)+u_0u_2u_4\sum\limits_{l=0}^{3s-1}\left( B_{2l}\prod\limits_{k_2=l+1}^{3s-1}A_{2k_2} \right)}{\left( \prod\limits_{k_1=0}^{3s}A_{2k_1}\right)+u_0u_2u_4\sum\limits_{l=0}^{3s}\left( B_{2l}\prod\limits_{k_2=l+1}^{3s}A_{2k_2} \right)},
\end{align}
\begin{align}
u_{6n+1}=&u_1 \frac{\left( \prod\limits_{k_1=0}^{3s-1}A_{2k_1+1}\right)+u_1u_3u_5\sum\limits_{l=0}^{3s-1}\left( B_{2l+1}\prod\limits_{k_2=l+1}^{3s-1}A_{2k_2+1} \right)}{\left( \prod\limits_{k_1=0}^{3s}A_{2k_1+1}\right)+u_1u_3u_5\sum\limits_{l=0}^{3s}\left( B_{2l+1}\prod\limits_{k_2=l+1}^{3s}A_{2k_2+1} \right)},
\end{align}
\begin{align}
u_{6n+2}=&u_2\prod_{s=0}^{n-1}\frac{\left( \prod\limits_{k_1=0}^{3s}A_{2k_1}\right)+u_0u_2u_4\sum\limits_{l=0}^{3s}\left( B_{2l}\prod\limits_{k_2=l+1}^{3s}A_{2k_2} \right)}{\left( \prod\limits_{k_1=0}^{3s+1}A_{2k_1}\right)+u_0u_2u_4\sum\limits_{l=0}^{3s+1}\left( B_{2l}\prod\limits_{k_2=l+1}^{3s+1}A_{2k_2} \right)},
\end{align}
\begin{align}
u_{6n+3}=&u_3 \prod_{s=0}^{n-1}\frac{\left( \prod\limits_{k_1=0}^{3s}A_{2k_1+1}\right)+u_1u_3u_5\sum\limits_{l=0}^{3s}\left( B_{2l+1}\prod\limits_{k_2=l+1}^{3s}A_{2k_2+1} \right)}{\left( \prod\limits_{k_1=0}^{3s+1}A_{2k_1+1}\right)+u_1u_3u_5\sum\limits_{l=0}^{3s+1}\left( B_{2l+1}\prod\limits_{k_2=l+1}^{3s+1}A_{2k_2+1} \right)},
\end{align}
\begin{align}
u_{6n+4}=&u_4 \prod_{s=0}^{n-1}\frac{\left( \prod\limits_{k_1=0}^{3s+1}A_{2k_1}\right)+u_0u_2u_4\sum\limits_{l=0}^{3s+1}\left( B_{2l}\prod\limits_{k_2=l+1}^{3s+1}A_{2k_2} \right)}{\left( \prod\limits_{k_1=0}^{3s+2}A_{2k_1}\right)+u_0u_2u_4\sum\limits_{l=0}^{3s+2}\left( B_{2l}\prod\limits_{k_2=l+1}^{3s+2}A_{2k_2} \right)},
\end{align}
\begin{align}
u_{6n+5}=&u_5 \prod_{s=0}^{n-1} \frac{\left( \prod\limits_{k_1=0}^{3s+1}A_{2k_1+1}\right)+u_1u_3u_5\sum\limits_{l=0}^{3s+1}\left( B_{2l+1}\prod\limits_{k_2=l+1}^{3s+1}A_{2k_2+1} \right)}{\left( \prod\limits_{k_1=0}^{3s+2}A_{2k_1+1}\right)+u_1u_3u_5\sum\limits_{l=0}^{3s+2}\left( B_{2l+1}\prod\limits_{k_2=l+1}^{3s+2}A_{2k_2+1} \right)},
\end{align}
whenever the denominators do not vanish, i.e.,
$u_0u_2u_4\sum_{l=0}^{s}B_{2l}\Theta^0(A_{k_2},l+1,s)\neq -\Theta^0(A_{k_1},0,s)$  
and
$u_1u_3u_5\sum_{l=0}^{s}B_{2l+1}\Theta^1(A_{k_2},l+1,s)\neq -\Theta^1(A_{k_1},0,s)$.
\end{subequations}
\par \noindent
The implication is that solutions of \eqref{1.1} are given by
\begin{subequations}\label{genesolxn}
\begin{align}
x_{6n-5}=&x_{-5} \prod_{s=0}^{n-1}\frac{\left( \prod\limits_{k_1=0}^{3s-1}a_{2k_1}\right)+x_{-5}x_{-3}x_{-1}\sum\limits_{l=0}^{3s-1}\left( b_{2l}\prod\limits_{k_2=l+1}^{3s-1}a_{2k_2} \right)}{\left( \prod\limits_{k_1=0}^{3s}a_{2k_1}\right)+x_{-5}x_{-3}x_{-1}\sum\limits_{l=0}^{3s}\left( b_{2l}\prod\limits_{k_2=l+1}^{3s}a_{2k_2} \right)},
\end{align}
\begin{align}
x_{6n-4}=&x_{-4} \prod_{s=0}^{n-1} \frac{\left( \prod\limits_{k_1=0}^{3s-1}a_{2k_1+1}\right)+x_{-4}x_{-2}x_{0}\sum\limits_{l=0}^{3s-
1}\left( b_{2l+1}\prod\limits_{k_2=l+1}^{3s-1}a_{2k_2+1} \right)}{\left( \prod\limits_{k_1=0}^{3s}a_{2k_1+1}\right)+x_{-4}x_{-2}x_{0}\sum\limits_{l=0}^{3s}\left( b_{2l+1}\prod\limits_{k_2=l+1}^{3s}a_{2k_2+1} \right)},
\end{align}
\begin{align}
x_{6n-3}=&x_{-3} \prod_{s=0}^{n-1}\frac{\left( \prod\limits_{k_1=0}^{3s}a_{2k_1}\right)+x_{-5}x_{-3}x_{-1}\sum\limits_{l=0}^{3s}\left( b_{2l}\prod\limits_{k_2=l+1}^{3s}a_{2k_2} \right)}{\left( \prod\limits_{k_1=0}^{3s+1}a_{2k_1}\right)+x_{-5}x_{-3}x_{-1}\sum\limits_{l=0}^{3s+1}\left( b_{2l}\prod\limits_{k_2=l+1}^{3s+1}a_{2k_2} \right)},
\end{align}
\begin{align}
x_{6n-2}=&x_{-2} \prod_{s=0}^{n-1}\frac{\left( \prod\limits_{k_1=0}^{3s}a_{2k_1+1}\right)+x_{-4}x_{-2}x_{0}\sum\limits_{l=0}^{3s}\left( b_{2l+1}\prod\limits_{k_2=l+1}^{3s}a_{2k_2+1} \right)}{\left( \prod\limits_{k_1=0}^{3s+1}a_{2k_1+1}\right)+x_{-4}x_{-2}x_{0}\sum\limits_{l=0}^{3s+1}\left( b_{2l+1}\prod\limits_{k_2=l+1}^{3s+1}a_{2k_2+1} \right)},
\end{align}
\begin{align}
x_{6n-1}=&x_{-1} \prod_{s=0}^{n-1}\frac{\left( \prod\limits_{k_1=0}^{3s+1}a_{2k_1}\right)+x_{-5}x_{-3}x_{-1}\sum\limits_{l=0}^{3s+1}\left( b_{2l}\prod\limits_{k_2=l+1}^{3s+1}a_{2k_2} \right)}{\left( \prod\limits_{k_1=0}^{3s+2}a_{2k_1}\right)+x_{-5}x_{-3}x_{-1}\sum\limits_{l=0}^{3s+2}\left( b_{2l}\prod\limits_{k_2=l+1}^{3s+2}a_{2k_2} \right)},
\end{align}
\begin{align}
x_{6n}=&x_0  \prod_{s=0}^{n-1} \frac{\left( \prod\limits_{k_1=0}^{3s+1}a_{2k_1+1}\right)+x_{-4}x_{-2}x_{0}\sum\limits_{l=0}^{3s+1}\left( b_{2l+1}\prod\limits_{k_2=l+1}^{3s+1}a_{2k_2+1} \right)}{\left( \prod\limits_{k_1=0}^{3s+2}a_{2k_1+1}\right)+x_{-4}x_{-2}x_{0}\sum\limits_{l=0}^{3s+2}\left( b_{2l+1}\prod\limits_{k_2=l+1}^{3s+2}a_{2k_2+1} \right)},
\end{align}
whenever the denominators do not vanish, i.e.,
$x_{-5}x_{-3}x_{-1}\sum_{l=0}^{s}b_{2l}\Theta^0(a_{k_2},l+1,s)\neq -\Theta^0(a_{k_1},0,s)$
and
$x_{-4}x_{-2}x_{-0}\sum_{l=0}^{s}b_{2l+1}\Theta^1(a_{k_2},l+1,s)\neq -\Theta^1(a_{k_1},0,s)$.
\end{subequations}
\subsection{The case where $(a_n)$ and $(b_n)$ are two-periodic sequences}
Let $a_n = (\lambda,\mu,\lambda,\mu,\lambda,\dots )$,  $b_n = (\eta,\zeta,\eta,\zeta,\dots )$, $ \Phi=x_{-5}x_{-3}x_{-1}$ and $ \Psi=x_{-4}x_{-2}x_{0}$. The solution in this case is given by the equations
\begin{subequations}
\begin{align}\label{genesolxn2p}
x_{6n-5}=&x_{-5} \prod_{s=0}^{n-1}\frac{\lambda^{3s}+\eta \Phi\sum\limits_{l=0}^{3s-1}\lambda^{l}}{\lambda^{3s+1}
+\eta \Phi\sum\limits_{l=0}^{3s}\lambda^{l}},
x_{6n-4}=x_{-4} \prod_{s=0}^{n-1} \frac{\mu^{3s}+\zeta \Psi\sum\limits_{l=0}^{3s-
1}\mu^l}{\mu^{3s+1}+\zeta \Psi\sum\limits_{l=0}^{3s}\mu^l},
\end{align}
\begin{align}
x_{6n-3}=&x_{-3} \prod_{s=0}^{n-1}\frac{\lambda^{3s+1}+\eta \Phi\sum\limits_{l=0}^{3s}\lambda^l}{\lambda^{3s+2}+
\eta \Phi\sum\limits_{l=0}^{3s+1}\lambda^l},
x_{6n-2}=x_{-2} \prod_{s=0}^{n-1}\frac{\mu^{3s+1}+\zeta \Psi\sum\limits_{l=0}^{3s}\mu^l}{\mu^{3s+2}+
\zeta \Psi\sum\limits_{l=0}^{3s+1}\mu^l},
\end{align}
\begin{align}
x_{6n-1}=&x_{-1} \prod_{s=0}^{n-1}\frac{\lambda^{3s+2}+\eta \Phi\sum\limits_{l=0}^{3s+1}\lambda^l}{\lambda^{3s+3}+\eta \Phi\sum\limits_{l=0}^{3s+2}\lambda^l },
x_{6n}=x_{0} \prod_{s=0}^{n-1} \frac{\mu^{3s+2}+\zeta \Psi\sum\limits_{l=0}^{3s+1}\mu^l}{\mu^{3s+3}+
\zeta \Psi\sum\limits_{l=0}^{3s+2}\mu^l},
\end{align}
\end{subequations}
where $\left(\sum\limits_{l=0}^{s-1}\lambda^l \right)\eta\Phi \neq -\lambda ^s $ and  $ \left(\sum\limits_{l=0}^{s-1}\lambda^l \right)\eta\Psi\neq -\lambda ^s,\quad s\leq 3n$.
\subsection{The case where $(a_n)$ and $(b_n)$ are constants}
Here, $a_n=\mu=\lambda$,  $b_n = \eta=\zeta$, $ \Phi=x_{-5}x_{-3}x_{-1}$ and $ \Psi=x_{-4}x_{-2}x_{0}$. Equation \eqref{1.1} becomes
$x_{n+1}={x_{n-5}}/{(\lambda +\eta x_{n-1}x_{n-3}x_{n-5})}$.
\subsubsection{$\lambda=1$}
The solution given in \eqref{genesolxn} simplifies to
\begin{subequations}\label{genesolxn1p}
\begin{align}
x_{6n-5}=&x_{-5} \prod_{s=0}^{n-1}\frac{1+(3s)\eta \Phi}{1
+(3s+1)\eta \Phi},
x_{6n-4}=x_{-4} \prod_{s=0}^{n-1} \frac{1+(3s)\eta \Psi}{1+(3s+1)\eta \Psi},
\end{align}
\begin{align}
x_{6n-3}=&x_{-3} \prod_{s=0}^{n-1}\frac{1+(3s+1)\eta \Phi}{1+
(3s+2)\eta \Phi},
x_{6n-2}=x_{-2} \prod_{s=0}^{n-1}\frac{1+
(3s+1)\eta \Psi}{1+
(3s+2)\eta \Psi},
\end{align}
\begin{align}
x_{6n-1}=&x_{-1} \prod_{s=0}^{n-1}\frac{1+{(3s+2)}\eta \Phi}{1+{(3s+3)}\eta \Phi },
x_{6n}=x_{0} \prod_{s=0}^{n-1} \frac{1+{(3s+2)}\eta \Psi}{1+{(3s+3)}\eta \Psi},
\end{align}
where $s\eta\Phi \neq -1 $ and  $ s\eta\Psi\neq -1,\quad s\leq 3n$.
\end{subequations}
\begin{itemize}
\item If  $\eta=1$, then equations in \eqref{genesolxn1p} are exactly the ones obtained in Theorem 2.1 in \cite{Elsayed}
 for
\begin{align}\label{1.1_11}
x_{n+1}=\frac{x_{n-5}}{1 + x_{n-1}x_{n-3}x_{n-5}}
\end{align}
 and their restriction (the initial conditions are arbitrary nonzero positive real numbers) is a special case of our restriction (the initial conditions are arbitrary nonzero real numbers and $s\Phi\neq -1,\,s\Psi \neq -1,\,0\leq s\leq 3n$ ).
\item
 If $ \eta=-1$, then equations in \eqref{genesolxn1p} are exactly the ones obtained in Theorem 3.1 in \cite{Elsayed} for
\begin{align}\label{1.1_1-1}
x_{n+1}=\frac{x_{n-5}}{1 - x_{n-1}x_{n-3}x_{n-5}}
\end{align}
and their restriction (the initial conditions are arbitrary nonzero real numbers and $jbdf\neq 1,\, jace\neq 1$ ) coincides with our restriction (the initial conditions are arbitrary nonzero real numbers and $s\Phi\neq 1,\,s\Psi \neq 1,\,0\leq s\leq 3n$ ).
\end{itemize}
\subsubsection{ The case where $\lambda\neq 1$}
In this case, the solution given in \eqref{genesolxn} simplifies to
\begin{align}\label{genesolxnnot1}
x_{6n-5}=&x_{-5} \prod_{s=0}^{n-1}\frac{\lambda^{3s}+\eta \Phi\left(\frac{1-\lambda ^{3s}}{1-\lambda}\right)}{\lambda^{3s+1}
+\eta \Phi\left(\frac{1-\lambda ^{3s+1}}{1-\lambda}\right)},
x_{6n-4}=x_{-4} \prod_{s=0}^{n-1} \frac{\lambda^{3s}+\zeta \Psi\left(\frac{1-\lambda ^{3s}}{1-\lambda}\right)}{\lambda^{3s+1}+\zeta \Psi\left(\frac{1-\lambda ^{3s+1}}{1-\lambda}\right)},\nonumber
\end{align}
\begin{align}
x_{6n-3}=&x_{-3} \prod_{s=0}^{n-1}\frac{\lambda^{3s+1}+\eta \Phi\left(\frac{1-\lambda ^{3s+1}}{1-\lambda}\right)}{\lambda^{3s+2}+
\eta \Phi\left(\frac{1-\lambda ^{3s+2}}{1-\lambda}\right)},
x_{6n-2}=x_{-2} \prod_{s=0}^{n-1}\frac{\lambda^{3s+1}+
\zeta \Psi\left(\frac{1-\lambda ^{3s+1}}{1-\lambda}\right)}{\lambda^{3s+2}+
\zeta \Psi\left(\frac{1-\lambda ^{3s+2}}{1-\lambda}\right)},\nonumber
\end{align}
\begin{align}
x_{6n-1}=&x_{-1} \prod_{s=0}^{n-1}\frac{\lambda^{3s+2}+\eta \Phi\left(\frac{1-\lambda ^{3s+2}}{1-\lambda}\right)}{\lambda^{3s+3}+\eta \Phi\left(\frac{1-\lambda ^{3s+3}}{1-\lambda}\right) },
x_{6n}=x_0 \prod_{s=0}^{n-1} \frac{\lambda^{3s+2}+\zeta \Psi\left(\frac{1-\lambda ^{3s+2}}{1-\lambda}\right)}{\lambda^{3s+3}+
\zeta \Psi\left(\frac{1-\lambda ^{3s+3}}{1-\lambda}\right)},
\end{align}
where $(1-\lambda ^s)\eta\Phi \neq -\lambda ^s(1-\lambda) $ and  $ (1-\lambda ^s)\zeta\Psi\neq- \lambda ^s(1-\lambda),\quad s\leq 3n$.\par \noindent
\textbf{Note.} If $\lambda=-1$ then  the
solution given in \eqref{genesolxnnot1} simplifies to
\begin{align}\label{genesolxn-1}
&x_{12n-5}=x_{6(2n)-5}=x_{-5},\;x_{12n-4}=x_{-4},\;x_{12n-3}=x_{-3},\;x_{12n-2}=x_{-2}, \nonumber\\&x_{12n-1}=x_{-1},
x_{12n}=x_{0},\; x_{12n+1}=
\frac{x_{-5}}{-1+\eta \Phi},\;
x_{12n+2}=
\frac{x_{-4}}{-1+\zeta\Psi},\;\nonumber\\&
 x_{12n+3}=
x_{-3}{(-1+\eta \Phi)},\;
x_{12n+4}=
x_{-2}{(-1+\zeta\Psi)},\; x_{12n+5}=
\frac{x_{-1}}{(-1+\eta \Phi)},\nonumber\\&
x_{12n+6}=
\frac{x_{0}}{-1+\zeta\Psi},\; x_{12n+7}=x_{6(2n+2)-5}
=x_{-5}.
\end{align}
\par \noindent
\begin{itemize}
\item
When setting $\eta =1$ in \eqref{genesolxn-1}, we get the result obtained in Theorem 4.1 in \cite{Elsayed} for
\begin{align}\label{1.1_-11'}
x_{n+1}=\frac{x_{n-5}}{-1 + x_{n-1}x_{n-3}x_{n-5}}
\end{align}
and their restriction coincides with our restriction  (the initial conditions are arbitrary nonzero real numbers, $\Phi \neq 1 $ and $ \Psi\neq1)$.
\item
When setting $\eta =-1$ in \eqref{genesolxn-1}, we get the result obtained in Theorem 5.1 in \cite{Elsayed} for
\begin{align}\label{1.1_-11}
x_{n+1}=\frac{x_{n-5}}{-1 - x_{n-1}x_{n-3}x_{n-5}}
\end{align}
and their restriction coincides with our restriction (the initial conditions are arbitrary nonzero real numbers, $\Phi \neq -1 $ and $ \Psi\neq -1)$.
\end{itemize}
\section{On the difference equations (\ref{un5'})}
Consider the fifth-order difference equations of the form (\ref{un5'}), i.e.,
\begin{align*}
u_{n+5}=\Omega=\frac{u_nu_{n+1}u_{n+2}}{u_{n+3}u_{n+4}(\lambda + \mu u_nu_{n+1}u_{n+2})}.
\end{align*}
Here, the procedure for finding the characteristics of \eqref{un5'} is similar as above and is as follows:
\begin{itemize}
\item [-]Impose the symmetry condition (\ref{LSC}) to \eqref{un5'}.
\item [-] Differentiate with respect to $u_n$ ( keeping $\Omega$ fixed) and viewing $u_{n+3}$ as a function of $u_n, u_{n+1},\; u_{n+2}$ and $\Omega.$
\item [-]Differentiate with respect to $u_n$ twice (keeping $u_{n+1}$ fixed).
\item [-]Use the method of separation.
\end{itemize}
After preforming this series of operations, we obtain the characteristics
\begin{align}\label{a6'}
\begin{split}
Q_5\left(n,{u_{n}}\right) =\beta^n  {u_n} \quad \text{and}\quad Q_6=\bar{\beta}^n  {u_n},
\end{split}
\end{align}
where $\beta= \exp{\left(\frac{-2\pi i}{3}\right)}$. Thus, we obtain two characteristics with corresponding generators given by
\begin{subequations}\label{gener}
\begin{eqnarray}
\begin{split}
X_5=&  \bar{\beta}^{n} u_n \partial u_n+  \bar{\beta}^{n+1} u_{n+1} \partial_{u_{ n+1}}+ \beta \bar{\beta}^{n} u_{n+2} \partial_{u_{ n+2}}+  \bar{\beta}^{n} u_{n+3} \partial_{u_{ n+3}}\\&+  \bar{\beta}^{n+1} u_{n+4} \partial_{u_{ n+4}},\label{gener1}
\end{split}
\\
\begin{split}
X_6=&  {\beta}^{n} u_n \partial u_n+ {\beta}^{n+1} u_{n+1} \partial_{u_{ n+1}}+ \bar{\beta}{\beta}^{n} u_{n+2} \partial_{u_{ n+2}}+ {\beta}^{n} u_{n+3} \partial_{u_{ n+3}}\\&+{\beta}^{n+1}u_{n+4} \partial_{u_{ n+4}}\label{gener2}.
\end{split}
\end{eqnarray}
\end{subequations}
From the characteristic equations 
\begin{align}\label{char1}
\dfrac{du_n}{{\bar{\beta}^n}u_n}=\frac{du_{n+1}}{\bar{\beta}^{n+1}u_{n+1}}
=\frac{du_{n+2}}{\bar{\beta}^{n+2}u_{n+2}}=\frac{du_{n+3}}{\bar{\beta}^{n+3}u_{n+3}}=
\frac{du_{n+4}}{\bar{\beta}^{n+4}u_{n+4}}\left(=\frac{V_n}{0}\right),
\end{align}
we obtain the invariants
$c_1= {u_{n+1}^{\beta}}/{u_n}, c_2= {u_{n+2}^{\bar{\beta}}}/{u_n}, c_3= u_nu_{n+3},  c_4={u_{n+4}^{\beta}}/{u_n}$  \text{and} $c_5=V_n.$
We readily notice that
\begin{align}\label{m1}
S^3\left( \dfrac{du_n}{{\bar{\beta}^n}u_n}=\frac{du_{n+1}}{\bar{\beta}^{n+1}u_{n+1}}\right)
=\frac{du_{n+3}}{\bar{\beta}^{n+3}u_{n+3}}=
\frac{du_{n+4}}{\bar{\beta}^{n+4}u_{n+4}}
\end{align}
and we choose
$V_n=c_1^{\bar{\beta}}c_2^{\beta}$,
i.e.,
\begin{align}\label{m4}
V_n={u_nu_{n+1}u_{n+2}}.
\end{align}
By shifting (\ref{m4}) thrice, we get
\begin{align}\label{m5}
V_{n+3}=\frac{V_n}{{\lambda}+{\mu}V_n}
\end{align}
whose solution is given by
\begin{equation}\label{m6}
V_n=
\begin{cases}
 \left(\frac{{\mu}\left[ ((-1)^{2/3}-1)n-(-1)^{1/3}+1+\beta\right]}{ 3\left((-1)^{2/3}-1\right)}+c_6 +\bar{\beta}^n c_7 +\beta ^n c_8\right)^{-1}\qquad \text{if}\qquad  {\lambda}= 1,\\ %\\
\left(c_{10}{\lambda}^{n/3}+\left[(-1)^{2/3}{\lambda}^{1/3}\right]^nc_8
+\left[-(-1)^{1/3}{\lambda}^{1/3}\right]^nc_9+ \frac{{\mu}}{1-{\lambda}} \right)^{-1} \text{ if } {\lambda}\neq 1.
\end{cases}
\end{equation}
The constants $c_i,\; i=6,\dots, 10$, can be obtained from the following equations:
\begin{subequations}
\begin{align}
&c_6+c_7+c_8=\frac{1}{u_0u_1u_2}-\frac{{\mu}\left[ -(-1)^{1/3}+1+\beta\right]}{ 3\left((-1)^{2/3}-1\right)}\label{567a}\\
&c_6+\bar{\beta}c_7+\beta c_8=\frac{1}{u_1u_2u_3}-\frac{{\mu}\left[ \beta-(-1)^{1/3}+(-1)^{2/3}\right]}{ 3\left((-1)^{2/3}-1\right)}\label{567b}
\\
&c_6+{\beta}c_7+\bar{\beta }c_8=\frac{1}{u_2u_3u_4}-\frac{{\mu}\left[ -1+2(-1)^{2/3}-(-1)^{1/3}+\beta\right]}{ 3\left((-1)^{2/3}-1\right)}\label{567c}\\
&c_{10}+c_8+c_9=\frac{1}{u_0u_1u_2}-\frac{{\mu}}{1-\lambda}\label{567d}\\
&\lambda ^{1/3}c_{10}+[(-1)^{2/3}\lambda ^{1/3}]c_8+[-(-1)^{1/3}\lambda ^{1/3}] c_9=\frac{1}{u_1u_2u_3}-\frac{{\mu}}{1-\lambda}\label{567e}\\
&\lambda ^{2/3}c_{10}+[(-1)^{2/3}\lambda ^{1/3}]^2c_8+[-(-1)^{1/3}\lambda ^{1/3}]^2 c_9=\frac{1}{u_2u_3u_4}-\frac{{\mu}}{1-\lambda}\label{567f}.
\end{align}
\end{subequations}
Thanks to (\ref{m4}), we can express $u_n$ in terms of $V_n$ as follows:
\begin{subequations}\label{m7}
\begin{align}
u_{n}=&\exp{\left( \beta ^n c_{11} + \bar{\beta}^n c_{12}% + \bar{\beta}^n c_2
-\frac{2}{\sqrt{3}}\left[\sum_{k=0}^{n-1}\text{Im}(\gamma (n,k))\ln V_{k}\right]\right)},\label{m9}
\end{align}
where $V_k$ is given in (\ref{m6}) with $\gamma(n,k)=\beta ^n \bar{\beta}^{k+1}$. The constants  $c_{11}$ and $c_{12}$ must satisfy
 \begin{align}\label{1112}
 c_{11}+c_{12}=\ln u_0,\quad
 \beta c_{11} +\bar{\beta} c_{12}=\ln u_1.
\end{align}
\end{subequations}
%\begin{rem}
Equations in (\ref{m7}) give the solutions of (\ref{un5'}) in a unified manner. \par \noindent For
%\end{rem}
 the sake of clarification, we split solutions (\ref{m9}) to realise the solutions in the existing literature. Using \eqref{m9} and (\ref{1112}) %together with the properties of $\gamma$ (given in \eqref{m11})
, we have
\begin{align}\label{m12}
\begin{split}
u_{6n}=&u_0 \prod _{s=1}^{2n}\frac{V_{3s-2}}{V_{3s-3}}.
\end{split}
\end{align}
Using the same approach, we have shown that
\begin{align}\label{m13}
u_{6n+i}=u_i \prod _{s=1}^{2n}\frac{V_{3(s-1)+i+1}}{V_{3(s-1)+i}},\quad i=0, \dots, 5.
\end{align}
\subsection{The case of ${\lambda}=1$}
Equation (\ref{un5'}) becomes
\begin{equation}\label{b1}
u_{n+5}=\frac{u_nu_{n+1}u_{n+2}}{u_{n+3}u_{n+4}(1+{\mu}u_nu_{n+1}u_{n+2})}.
\end{equation}
and we said earlier that the solution of (\ref{m5}), in this case, is (\ref{m6}), i.e.,
\begin{align}\label{m6a}
V_n= \left(\frac{{\mu}\left[ ((-1)^{2/3}-1)n-(-1)^{1/3}+1+\beta\right]}{ 3\left((-1)^{2/3}-1\right)}+c_6 +\bar{\beta}^n c_7 +\beta ^n c_8\right)^{-1}.
\end{align}
We have 
\begin{align}\label{m14a}
V_{3s}=&\left({\frac{{\mu}\left[ ((-1)^{2/3}-1)(3s)-(-1)^{1/3}+1+\beta\right]}{ 3\left((-1)^{2/3}-1\right)}+c_6 + c_7 + c_8}\right)^{-1}%.
\end{align}
and using (\ref{567a}) in (\ref{m14a}), we get
\begin{align}\label{m14bb}
V_{3s}=&
\frac{u_0u_1u_2}{1+{\mu}su_0u_1u_2}.
\end{align}
Using the same approach, we have shown that
\begin{flushleft}
\begin{align}\label{a=1}
\begin{split}
V_{3s}=\frac{u_0u_1u_2}{1+{\mu}su_0u_1u_2};\;
V_{3s+1}=\frac{u_1u_2u_3}{1+{\mu}su_1u_2u_3};\;
V_{3s+2}=\frac{u_2u_3u_4}{1+{\mu}su_2u_3u_4}.
\end{split}
\end{align}
Let $a=x_{-4},\;b=x_{-3},\;c=x_{-2},\;d=x_{-1}$ and $A=u_0,\;B=u_1$, $C=u_2,$ $D=u_3$ and $E=u_4$. {Using (\ref{a=1}) in (\ref{m13}), we obtain the solution of (\ref{b1}) as follows:}
\begin{align}\label{a=11}
u_{6n}=&\frac{D^{2n}}{A^{2n-1}} \prod _{s=1}^{2n-1}\frac{1+{\mu}sABC}{1+{\mu}sBCD}, 
u_{6n+1}=\frac{E^{2n}}{B^{2n-1}} \prod _{s=1}^{2n-1}\frac{1+{\mu}sBCD}{1+{\mu}sCDE},\nonumber\\ 
u_{6n+2}=&\frac{CA^{2n}B^{2n}}{D^{2n}E^{2n}} \prod _{s=0}^{2n-1}\frac{1+{\mu}sCDE}{1+{\mu}(s+1)ABC}, 
u_{6n+3}=\frac{D^{2n+1}}{A^{2n}} \prod _{s=1}^{2n}\frac{1+{\mu}sABC}{1+{\mu}sBCD} ,\nonumber\\ 
u_{6n+4}=&\frac{E^{2n+1}}{B^{2n}} \prod _{s=0}^{2n}\frac{1+{\mu}sBCD}{1+{\mu}sCDE}, 
u_{6n+5}=\frac{C(AB)^{2n+1}}{(DE)^{2n+1}} \prod _{s=0}^{2n}\frac{1+{\mu}sCDE}{1+{\mu}(s+1)ABC}
\end{align}
whenever the denominators do not vanish.
\end{flushleft}
\subsection{The case of ${\lambda} \neq 1$}
In this case, as we found earlier, the solution of (\ref{m5}) is given by  (\ref{m6}), i.e.,
\begin{align}\label{m6aa}
V_n=\left( c_{10} {\lambda}^{n/3}+[(-1)^{2/3}{\lambda}^{1/3}]^nc_8+[-(-1)^{1/3}{\lambda}^{1/3}]^nc_9+
\frac{{\mu}}{1-{\lambda}}\right)^{-1}.
\end{align}
Using this, we get
\begin{align}\label{m14aa}
V_{3s}=&\left(c_{10}{\lambda}^s + c_8{\lambda}^s + c_9{\lambda}^s+\frac{{\mu}}{1-{\lambda}}\right)^{-1}
\end{align}
and using (\ref{567d}) in (\ref{m14aa}), we find
\begin{align}\label{m14bba}
V_{3s}=&
\frac{ABC}{{\lambda}^s+{\mu}ABC\left( \frac{1-{\lambda}^s}{1-{\lambda}}\right)}.
\end{align}
{Using the same approach, we have shown that }
\begin{align}\label{a=1a}
V_{3s+i}=&\frac{u_iu_{i+1}u_{i+2}}{{\lambda}^s+{\mu}u_iu_{i+1}u_{i+2}\left( \frac{1-{\lambda}^s}{1-{\lambda}}\right)},\quad i=0,1,2.
\end{align}
{Using (\ref{a=1a}) in (\ref{m13}), we obtain the solution of (\ref{un5'}) as follows:}
\begin{align}\label{a=11a'}
u_{6n}=\frac{D^{2n}}{A^{2n-1}} \prod _{s=2}^{2n}\frac{{\lambda}^{s-1}+{\mu}\Delta_0\sum \limits_{j=0}^{s-2}\lambda ^j}{{\lambda}^{s-1}+{\mu}\Delta_1 \sum\limits_{j=0}^{s-2}\lambda ^j},
u_{6n+1}=\frac{E^{2n}}{B^{2n-1}} \prod _{s=2}^{2n}\frac{{\lambda}^{s-1}+{\mu}\Delta_1 \sum\limits_{j=0}^{s-2}\lambda ^j}{{\lambda}^{s-1}+{\mu}\Delta_2 \sum\limits_{j=0}^{s-2}\lambda ^j}, \nonumber
\end{align}
\begin{align}
u_{6n+2}=C\frac{\Delta_0^{2n}}{\Delta_1^{2n}} \frac{\prod \limits _{s=2}^{2n}\left({\lambda}^{s-1}+{\mu}\Delta_2 \sum\limits_{j=0}^{s-2}\lambda ^j\right)}{\prod \limits_{s=1}^{2n}\left({\lambda}^s+{\mu}\Delta_0 \sum\limits_{j=0}^{s-1}\lambda ^j\right)} ,
u_{6n+3}=&\frac{D^{2n+1}}{A^{2n}} \prod _{s=2}^{2n+1}\frac{{\lambda}^{s-1}+{\mu}\Delta_0 \sum\limits_{j=0}^{s-2}\lambda ^j}{{\lambda}^{s-1}+{\mu}\Delta_1 \sum\limits_{j=0}^{s-2}\lambda ^j},\nonumber
\end{align}
\begin{align}
u_{6n+4}=&E\frac{\Delta_2^{2n}}{\Delta_1^{2n}} \prod\limits _{s=2}^{2n+1}\frac{{\lambda}^{s-1}+{\mu}\Delta_1 \sum\limits_{j=0}^{s-2}\lambda ^j}{{\lambda}^{s-1}+{\mu}\Delta_2\sum\limits_{j=0}^{s-2}\lambda ^j},  
u_{6n+5}=C\frac{\Delta_0^{2n+1}}{\Delta_1^{2n+1}}\frac{\prod \limits _{s=2}^{2n+1}\left({\lambda}^{s-1}+{\mu}\Delta_2 \sum\limits_{j=0}^{s-2}\lambda ^j\right)}{\prod\limits _{s=1}^{2n+1}\left({\lambda}^s+{\mu}\Delta_0 \sum\limits_{j=0}^{s-1}\lambda ^j\right)},
\end{align}
where $\Delta_i=u_iu_{i+1}u_{i+2}$.
Equations in (\ref{a=11a'}) give the exact solution of (\ref{un5'}) for any real values of $\lambda$ and $\mu$ provided that the denominators do not vanish.
\par \noindent
Recall that we acted the shift operator on (\ref{xn5}) to get (\ref{un5'}). Hence, the solutions of (\ref{xn5}) are obtained, using (\ref{a=11a'}), as follows:
\begin{align}\label{Xa=11a'}
x_{6n-4}=&\frac{d^{2n}}{a^{2n-1}} \prod _{s=2}^{2n}\frac{{\lambda}^{s-1}+{\mu}abc\sum\limits_{j=0}^{s-2}\lambda ^j}{{\lambda}^{s-1}+{\mu}bcd \sum\limits_{j=0}^{s-2}\lambda ^j}, 
x_{6n-3}=\frac{e^{2n}}{b^{2n-1}} \prod _{s=2}^{2n}\frac{{\lambda}^{s-1}+{\mu}bcd \sum\limits_{j=0}^{s-2}\lambda ^j}{{\lambda}^{s-1}+{\mu}cde \sum\limits_{j=0}^{s-2}\lambda ^j},\nonumber
\end{align}
\begin{align}
x_{6n-2}=&\frac{c(ab)^{2n}}{(de)^{2n}} \frac{\prod \limits _{s=2}^{2n}\left({\lambda}^{s-1}+{\mu}cde \sum\limits_{j=0}^{s-2}\lambda ^j\right)}{\prod\limits _{s=1}^{2n}\left({\lambda}^s+{\mu}abc \sum\limits_{j=0}^{s-1}\lambda ^j\right)},
x_{6n-1}=\frac{d^{2n+1}}{a^{2n}} \prod\limits _{s=2}^{2n+1}\frac{{\lambda}^{s-1}+{\mu}abc \sum\limits_{j=0}^{s-2}\lambda ^j}{{\lambda}^{s-1}+{\mu}bcd \sum\limits_{j=0}^{s-2}\lambda ^j}\nonumber , 
\end{align}
\begin{align}
x_{6n}=&\frac{e^{2n+1}}{b^{2n}} \prod _{s=2}^{2n+1}\frac{{\lambda}^{s-1}+{\mu}bcd \sum\limits_{j=0}^{s-2}\lambda ^j}{{\lambda}^{s-1}+{\mu}cde\sum\limits_{j=0}^{s-2}\lambda ^j}, 
x_{6n+1}=\frac{c(ab)^{2n+1}}{(de)^{2n+1}}\frac{\prod\limits _{s=2}^{2n+1}\left[{\lambda}^{s-1}+{\mu}cde \sum\limits_{j=0}^{s-2}\lambda ^j\right]}{\prod _{s=1}^{2n+1}\left({\lambda}^s+{\mu}abc \sum\limits_{j=0}^{s-1}\lambda ^j\right)},
\end{align}
 for any real values of $\lambda$ and $\mu$ as long as the denominators do not vanish.
\begin{itemize}
\item
 When $\lambda=1$ and ${\mu}=1$, equations in (\ref{Xa=11a'}) yield the results obtained by Yasin Yazlik in Theorem 5 in \cite{yy} for
\begin{align}\label{xna1b1'}
x_{n+1}=\frac{x_{n-2}x_{n-3}x_{n-4}}{x_nx_{n-1}(1+x_{n-2}x_{n-3}x_{n-4})},\quad n=0,1,2,\dots
\end{align}
where $a,\; b,\;c,\;c,\; d$ and $e$ are positive real numbers.
\item
When $\lambda =1$ and ${\mu}=-1$, equations in (\ref{Xa=11a'}) yield the results obtained by Yasin Yazlik in Theorem 9 in \cite{yy} for
\begin{align}\label{xna1b1}
x_{n+1}=\frac{x_{n-2}x_{n-3}x_{n-4}}{x_nx_{n-1}(1-x_{n-2}x_{n-3}x_{n-4})},\quad n=0,1,2,\dots
\end{align}
where $a,\; b,\;c,\;c,\;d$ and $e$ are positive real numbers with ${a}{b}{c}\neq 1$   and  ${c}{d}{e}\neq 1$ .\par \noindent
\textbf{Note. }There should not be a minus sign right after the expression of $x_{3n-2}$ in Theorem 9 in \cite{yy}.
\item
When $\lambda=-1$ and $\mu=1$, equations in (\ref{Xa=11a'}) yield the results obtained by Yasin Yazlik in Theorem 7 in \cite{yy} for
\begin{align}\label{xna1b1a'}
x_{n+1}=\frac{x_{n-2}x_{n-3}x_{n-4}}{x_nx_{n-1}(-1+x_{n-2}x_{n-3}x_{n-4})},\quad n=0,1,2,\dots
\end{align}
where $a,\; b,\;c,\; d$ and $e$ are non zero real numbers with $a bc\neq 1$,\; $bcd\neq 1$ and $c de\neq1$.
\par \noindent
\item
When $\lambda=-1$ and $\mu=-1$, equations in (\ref{Xa=11a'}) yield the results obtained by Yasin Yazlik in Theorem 11 in \cite{yy} for
\begin{align}\label{xna1b1a}
x_{n+1}=\frac{x_{n-2}x_{n-3}x_{n-4}}{x_nx_{n-1}(-1-x_{n-2}x_{n-3}x_{n-4})},\quad n=0,1,2,\dots
\end{align}
where $a,\; b,\;c,\;c,\;d$ and $e$ are non zero real numbers with ${a}{b}{c}\neq - 1$ ,\;${b}{c}{d}\neq - 1$  and  ${c}{d}{e}\neq -1$ .\par \noindent
\textbf{Note. }There should be a minus sign right after the expression of $x_{6n+1}$ in Theorem 11 in \cite{yy}.
\end{itemize}
\section{Conclusion}
In this paper, we have obtained nontrivial symmetries of some rational ordinary difference equations and their exact solutions were obtained. Most importantly, we note that \eqref{ab} gives a clear idea, without making any lucky guesses, for the most convenient choice of the invariant of \eqref{1.2}.
\section*{Conflict of interests}
The authors declare that there is no conflict of interest regarding the publication of this paper.
\section*{Data availability statement}
No data were used to support this study.
\section*{Funding statement}
The authors received no specific funding for this work.

\end{document}